\documentclass{article}

\usepackage{amsmath}
\usepackage{amsfonts}
\usepackage{amscd}
\usepackage{amsthm}
\usepackage{graphics}
\usepackage{verbatim}
\usepackage{enumerate}
\newtheorem{thm}{Theorem}[subsection]
\newtheorem{prop}[thm]{Proposition}
\newtheorem{lem}[thm]{Lemma}

\theoremstyle{definition}

\newtheorem{rem}[thm]{Remark}

\newcommand{\C}{\mathbb{C}}
\newcommand{\Z}{\mathbb{Z}}
\newcommand{\Mabc}{M_{(a, b, c)}}
\newcommand{\Fabc}{\Phi_{(a, b, c)}}
\newcommand{\Md}{M_{(0, 2, 2)}}
\newcommand{\Mdt}{\widetilde{M}_{(0, 2, 2)}}
\newcommand{\Q}{\mathbb{Q}}
\newcommand{\St}{\widetilde{\Sigma}}
\newcommand{\M}{\mathbf{M}}
\newcommand{\Muo}{\widetilde{\mu}(\omega+\overline{\omega})}
\newcommand{\om}{(\omega+\overline{\omega})}

\begin{document}
\title{
       On the Beauville form of the known \\
       irreducible symplectic varieties}
\author{Antonio Rapagnetta\\%\thanks{support}\\
        Dipartimento di Matematica Universit\`a di Pisa\\
        rapagnetta@mail.dm.unipi.it}

\date{May 15, 2006}
\maketitle

\begin{abstract}
We study the global geometry of the ten dimensional O'Grady irreducible symplectic variety.
We determine its second Betti number, its Beauville form and its Fujiki constant.

\end{abstract}
\section*{Introduction}
Irreducible symplectic varieties are simply connected compact
K\"ahler manifolds with a unique, up to $\mathbb{C}^{*}$, global
holomorphic two form and such that this two form is nondegenerate
at each point. By Bogomolov's decomposition theorem \cite{Bo74},
they are fundamental   in the classification of compact K\"ahler
manifolds with torsion $c_{1}$.

Very few examples of irreducible
symplectic varieties are known. For any positive integer $n$,
Beauville exhibited two examples of dimension $2n$ (\cite{Be83}):
the Hilbert scheme $X^{[n]}$, parametrizing $0$-dimensional
subschemes of length $n$ on a K3 surface $X$, and the generalized
Kummer variety ${\rm K}^{n}(T)$ of a $2$-dimensional torus $T$.
More recently O'Grady,  in \cite{OG99} and \cite{OG03}, found two
new examples $\widetilde{\mathcal{M}}$ and $\mathbf{M}$ of
dimension six and ten. Up to now, each known  example   can be
deformed into one of these examples.

For any irreducible symplectic variety $Y$, the group
$H^{2}(Y,\Z)$ is endowed with a deformation invariant  integral
primitive bilinear form $B_{Y}$, called the Beauville form.
Related with the Beauville form there is a positive rational
constant $c_{Y}$, called the Fujiki constant, which is a
topological invariant of $Y$. The Beauville form and the Fujiki
constant are fundamental invariants of an irreducible symplectic
variety: in the cases of the Beauville examples they have been
known since the paper \cite{Be83} and in the case of the six
dimensional O'Grady example $\widetilde{\mathcal{M}}$ they were
determined in \cite{Ra04}.

In this paper we deal with the remaining case:
we determine the second Betti number of the
ten dimensional O'Grady example $\mathbf{M}$ and then compute
$B_{\mathbf{M}}$ and $c_{\mathbf{M}}$.
   In the following table we give the complete
list of the Beauville forms and the Fujiki constants
of all known irreducible symplectic verieties.

\[
\begin{array}{|c|c|c|c|c|}
\hline
Y & {\rm dim}(Y) & b_{2}(Y) & c_{Y} & (H^{2}(Y,\Z), B_{Y}) \\
\hline
X^{[n]} & 2n & 23 & \frac{(2n)!}{n!2^{n}}& H^{\oplus_{\perp} 3}\oplus_{\perp} -E_{8}^{\oplus_{\perp} 3}\oplus_{\perp}(-2(n-1)) \\
\hline
{\rm K}^{n}(T) & 2n & 7 & \frac{(2n)!}{n!2^{n}}(n+1)& H^{\oplus_{\perp} 3}
\oplus_{\perp}(-2(n+1)) \\
\hline
\widetilde{\mathcal{M}}& 6 & 8 & 60& H^{\oplus_{\perp} 3}
\oplus_{\perp}(-2)^{\oplus_{\perp} 2} \\
\hline
\mathbf{M} & 10 & 24 & 945 & H^{\oplus_{\perp} 3}\oplus_{\perp} -E_{8}^{\oplus_{\perp} 3}\oplus_{\perp} \Lambda \\
\hline
\end{array}
\]
In this table the lattice $H$ is the
standard hyperbolic plane, the lattice $-E_{8}$
is the unique negative definite even unimodular
lattice of rank eight and $(i)$ is the rank $1$
lattice generated by an element whose square is $i$.
Finally $\Lambda$ is a rank $2$ lattice whose
associated matrix in a suitable basis is
$$\left(\begin{matrix}
-6 & 3 \\
3 & -2  \\
\end{matrix}\right).$$

It is remarkable that both the O'Grady examples have the same
Fujiki constants of Beauville examples of the same dimension: the
Fujiki constant of $\widetilde{\mathcal{M}}$ equals the Fujiki
constant of ${\rm K}^{3}(T)$ and the Fujiki constant of  $\M$
equals the one of $X^{[5]}$.

{\bf Acknowledgements.} I would like to thank  Kieran O'Grady
for useful conversations and Donatella Iacono and Francesco Esposito
for their helpful support.

\section{$b_{2}(\M)=24$}
In this paper $X$ is a K3 surface such that $Pic(X)=\Z<H>$ and
$H^{2}=2$. We denote by $\Mabc$ the Simpson moduli space of
semistable sheaves on $X$ with  Mukai vector
$(a,bc_{1}(H),c\eta)\in H^{0}(X,\Z)\oplus H^{2}(X,\Z)\oplus
H^{4}(X,\Z)$, where $\eta$ is the fundamental form in
$H^{4}(X,\Z)$.

If $a=0$ and $b>0$ the moduli space $\Mabc$ parametrizes pure
1-dimensional sheaves on $X$, hence there exists a regular
morphism $\Fabc:\Mabc\rightarrow|bH|$ sending a sheaf to its
Fitting subscheme (see \cite{LP93} end \cite{Ei95}).
The general fiber of $\Fabc$ is
easily described: if $C\in|bH|$ is smooth then $\Fabc^{-1}(C)$
parametrizes push-forwards in $X$ of line bundles of degree $b^{2}+c$ on $C$.
More
precisely there is an isomorphism $\Fabc^{-1}(C)\simeq
Pic^{b^{2}+c}(C)$.

If $a=0$, $b=2$ and $c$ is odd, the moduli space $\Mabc$ is an
irreducible symplectic variety which can be deformed into the
Hilbert scheme $X^{[5]}$  parametrizing 0-dimensional subschemes
of length $5$ on $X$ (see \cite{Yo}) and
$\Fabc:\Mabc\rightarrow|2H|$ is a Lagrangian fibration by
Matsushita's theorem (see \cite{Ma01}). On the other hand $\Md$ is
an irreducible variety (see Theorem 4.4 of \cite{KLS05}) that admits
a symplectic resolution $\widetilde{\pi}:\Mdt\rightarrow\Md$ and
is  birational to the 10-dimensional O'Grady example $\mathbf{M}$
(see Proposition (4.1.5) of \cite{OG99} and Remark 1.1.9 of \cite{Ra04}).
The composition
$\Psi:=\Phi_{(0,2,2)}\circ\widetilde{\pi}:\Mdt \rightarrow|2H|$ gives also
in this case a Lagrangian fibration on the irreducible symplectic
variety $\Mdt$.

By a theorem of  Huybrechts (see Theorem 4.6 of \cite{Hu99}) $\Mdt$
and $\mathbf{M}$ are deformation equivalent, so they have the same
Betti numbers. We are going to determine $b_{2}(\Mdt)$ by
comparing the Lagrangian fibrations $\Phi:=\Phi_{(0,2,1)}$ and
$\Psi$ and then using results about the topology of $M_{(0,2,1)}$.
\begin{thm}\label{bdue}
The dimension of $H^{2}(\M,\Q)$ is $24$.
\end{thm}
This theorem is a consequence of the following three propositions.
The first proposition compares open subsets of  $\Mdt$ and $M_{(0,2,1)}$.
\begin{prop}\label{confronto aperti}
Let $U\subset|2H|$ be the open subset parametrizing irreducible curves.
%curves with at
%most nodal singularities and singular in at most a unique point.
The vector spaces  $H^{2}(\Psi^{-1}(U),\Q)$ and
$H^{2}(\Phi^{-1}(U),\Q)$ are isomorphic.
\end{prop}
In the second proposition we determine the second Betti number of
the open subvariety $\Phi^{-1}(U)\subset M_{(0,2,1)}$.
\begin{prop}\label{coomologia aperto}
The dimension of $H^{2}(\Phi^{-1}(U),\Q)$ is $22$.
\end{prop}
In the third proposition we collect the required informations on the
complement of $\Psi^{-1}(U)$ in $\mathbf{M}$.
\begin{prop}\label{componenti irriducibili}
Let $R\subset|2H|$ be the locus parametrizing reducible curves.
The divisor  $\Psi^{-1}(R)\subset\M$ is the union of two
irreducible divisors.
\end{prop}
In the remaining part of this section, coefficients of singular
cohomology groups are always rational and, for simplicity,
we often omit them in the notation.

Before proving the three propositions we assume them to prove
Theorem \ref{bdue}.

\noindent{\em Proof of Theorem \ref{bdue}.} Let $Y$ be the
singular locus of $\Psi^{-1}(R)$, the couple $(\Mdt\setminus
Y, \Psi^{-1}(U))$ induces the following  exact sequence of
cohomology groups with rational coefficients
$$H^{2}(\Mdt\setminus Y, \Psi^{-1}(U)) \rightarrow
H^{2}(\Mdt\setminus Y) \rightarrow H^{2}(\Psi^{-1}(U)).$$
Let $N$ be the normal bundle of $\Psi^{-1}(R)\setminus Y$ in
$\Mdt\setminus Y$ and let $N^{0}$ be the complement of its zero
section: by excision $H^{2}(\Mdt\setminus Y,
\Psi^{-1}(U))= H^{2}(N,N^{0})$, by Thom isomorphism
$H^{2}(N,N^{0})=H^{0}(\Psi^{-1}(R)\setminus Y)$ and by Proposition
\ref{componenti irriducibili} we get $H^{0}(\Psi^{-1}(R)\setminus
Y)=\Q^{2}$.

Since $H^{2}(\Psi^{-1}(U))=H^{2}(\Phi^{-1}(U))=\Q^{22}$ by
Proposition \ref{confronto aperti} and Proposition \ref{coomologia
aperto}, the previous exact sequence implies $b_{2}(\Mdt\setminus
Y)\le 24$. Since  $Y$ has codimension two in $\Mdt$ we also get
$b_{2}(\Mdt)\le 24$. Finally  $b_{2}(\M)=24$ since $b_{2}(\M)=b_{2}(\Mdt)$ and
O'Grady
already proved $b_{2}(\mathbf{M})\ge 24$ (see
\cite{OG99}).\qed \medskip

\begin{rem} \label{rivestimento doppio}
Under our assumption the linear system $|H|$ induces a double covering
$f:X\rightarrow \mathbb{P}^{2}(\simeq |H|^{\vee})$ ramified over a
smooth sextic. In what follows we often use that $f^{*}$ induces a
bijection between plane conics and curves belonging to the linear
system $|2H|$.
\end{rem}
\noindent{\em Proof of Proposition \ref{confronto aperti}.} Denote
by $U_{0}\subset U$ the open subset parametrizing curves whose
singular loci are empty or consist of a unique nodal point. Since
the complement of $U_{0}$ in $U$ has codimension two and the
fibers of $\Psi$ and $\Phi$ are equidimensional (see \cite{Ma00}),
the same property holds for the complements of $\Psi^{-1}(U_{0})$
and $\Phi^{-1}(U_{0})$ in $\Psi^{-1}(U)$ and $\Phi^{-1}(U)$
respectively. Hence Proposition \ref{confronto aperti} is
equivalent to the equality
$b_{2}(\Psi^{-1}(U_{0}))=b_{2}(\Phi^{-1}(U_{0}))$.

Let $\Psi_{0}:\Psi^{-1}(U_{0})\rightarrow U_{0}$ and $\Phi_{0}:
\Phi^{-1}(U_{0})\rightarrow U_{0}$ be the restrictions of $\Psi$ and
$\Phi$.
The $E_{2}^{p,q}$ terms of the associated Leray spectral sequences are
$H^{p}(R^{q}\Psi_{0*}(\Q))$ and
$H^{p}(R^{q}\Phi_{0*}(\Q))$.

Since the abutment of the Leray spectral sequence is the cohomology
of the domain, Proposition \ref{confronto aperti} follows if we prove:
\begin{enumerate}[a)]
\item{$R^{q}\Psi_{0*}(\Q)$ and $R^{q}\Phi_{0*}(\Q)$ are isomorphic
for $q\le2$.}
\item{For $p+q=2$, the $E_{2}^{p,q}$ terms of both the spectral
sequences survive to
the $E_{\infty}$ pages.}
\end{enumerate}

Let's prove a).
Let $U_{s}\subset |2H|$ be the locus parametrizing smooth curves and
let $i:U_{s}\rightarrow U_{0}$ be the open inclusion.
Statement a) is an obvious consequence of
\begin{enumerate}
\item{For $q\le 2$ the sheaves $i^{*}R^{q}\Psi_{0*}(\Q)$ and $i^{*}R^{q}\Phi_{0*}(\Q))$ are isomorphic.}
\item{For $q\le 2$  the natural attachment  maps
$R^{q}\Psi_{0*}(\Q)\rightarrow i_{*}i^{*}R^{q}\Psi_{0*}(\Q)$ and
$R^{q}\Phi_{0*}(\Q)\rightarrow i_{*}i^{*}R^{q}\Phi_{0*}(\Q)$ are
isomorphisms.}
\end{enumerate}
In order to prove 1, we denote by $\Psi_{s}:\Psi^{-1}(U_{s})\rightarrow
U_{s}$ and $\Phi_{s}:
\Phi^{-1}(U_{s})\rightarrow U_{s}$ the restrictions of $\Psi$ and
$\Phi$.
Since $R^{q}\Psi_{s*}(\Q)\simeq i^{*}R^{q}\Psi_{0*}(\Q)$
and $R^{q}\Phi_{s*}(\Q)\simeq i^{*}R^{q}\Phi_{0*}(\Q)$
we are reduced to show that there exists an isomorphism between
$R^{q}\Psi_{s*}(\Q)$ and $R^{q}\Phi_{s*}(\Q)$.

Let $l\in |2H|\simeq \mathbb{P}^{5}$ be a general line, by the Zariski
theorem the inclusion induces a surjection on fundamental groups
$\pi_{1}(l\cap U_{s})\rightarrow \pi_{1}(U_{s})$.
Since $\Psi_{s}$ and $\Phi_{s}$ are smooth
$R^{q}\Psi_{s*}(\Q)$ and $R^{q}\Phi_{s*}(\Q)$
are local systems, hence they are isomorphic if and only if their
restrictions $R^{q}\Psi_{s*}(\Q)_{|l\cap U_{s}}$ and
$R^{q}\Phi_{s*}(\Q)_{|l\cap U_{s}}$
to $l\cap U_{s}$ are isomorphic.
Hence, denoting by $\Psi_{l}:\Psi^{-1}(l\cap U_{s})\rightarrow l\cap U_{s}$
and $\Phi_{l}:
\Phi^{-1}(l\cap U_{s})\rightarrow l\cap U_{s}$
the restrictions of $\Psi$ and $\Phi$ we want an isomorphism between
$R^{q}\Psi_{l*}(\Q)$ and $R^{q}\Phi_{l*}(\Q)$.
Since $l\cap U_{s}$ parametrizes smooth curves
the family $\Psi_{l}:\Psi^{-1}(l\cap U_{s})\rightarrow l\cap U_{s}$ is
isomorphic to the degree $6$ relative Picard group ${\rm Pic}^{6}(l\cap
U_{s})\rightarrow l\cap U_{s}$ of the family of curves parametrized
by $l\cap U_{s}$ and analogously $\phi_{l}$ can be identified with the
degree $5$ relative Picard group ${\rm Pic}^{5}(l\cap
U_{s})\rightarrow l\cap U_{s}$ of the same family of curves.
The wanted isomorphism exists since ${\rm Pic}^{6}(l\cap
U_{s})\simeq{\rm Pic}^{5}(l\cap
U_{s})$ over $l\cap U_{s}$. This follows since the family of
curves parametrized by $l\cap U_{s}$ admits sections: any point
in the base locus of $l$ gives a section.

%\sum_{p+q=2}h^{p}(R^{q}\overline{\Phi}_{*}(\Q))=
%\sum_{p+q=2}h^{p}(R^{q}\overline{\Psi}_{*}(\Q))=
%h^{2}(\Phi^{-1}(U),\Q).\]
In order to prove 2 we need a general Lemma.
\begin{lem}\label{lemma degenerazione}
Let $g:\mathcal{X}\rightarrow \Delta$ be a proper map with
irreducible fibers, from a complex smooth surface $\mathcal{X}$
onto the open unit disk $\Delta\subset\C$. Suppose that $g$ has a
unique critical point $p$, suppose that $p$ is non degenerate and
$g(p)=0$. Let
$\hat{g}:\overline{Pic}^{i}(\mathcal{X})\rightarrow\Delta$ be the
compactification, by torsion free sheaves, of the degree $i$
relative Picard group $Pic^{i}(\mathcal{X})$. Let
$i:\Delta^{*}:=\Delta \setminus {0}\rightarrow\Delta$ be the
inclusion. For $q\le 2$, the natural attachment map
$\alpha:R^{q}\hat{g}_{*}\Q\rightarrow
i_{*}i^{*}R^{q}\hat{g}_{*}\Q$ is an isomorphism.
\end{lem}
\begin{proof}
We only have to check that the map $\Gamma(\alpha):
\Gamma(R^{q}\hat{g}_{*}\Q)\rightarrow
\Gamma(i_{*}i^{*}R^{q}\hat{g}_{*}\Q)$
induced by $\alpha$ on global
sections is an isomorphism.

The vector space $\Gamma(R^{q}\hat{g}_{*}\Q)$ is  isomorphic to $
H^{q}(\overline{Pic}^{i}(\mathcal{X}))$ and the vector space
$\Gamma(i_{*}i^{*}R^{q}\hat{g}_{*}\Q)$ is isomorphic to the group
$H^{q}(\hat{g}^{-1}(\frac{1}{2}))^{\pi_{1}(\Delta^{*})}$ of the
$q$-cocycles of the general fiber of $\hat{g}$ that are invariant under the
monodromy action of $\pi_{1}(\Delta^{*})$. Finally, using these
identifications, the map $\Gamma(\alpha)$ is just the map induced
in cohomology by the inclusion
$\hat{g}^{-1}(\frac{1}{2})\hookrightarrow
\overline{Pic}^{i}(\mathcal{X})$.

The central fiber of $\hat{g}$ is a normal crossings divisor (see
\cite{Se00}) and $\hat{g}$ is a semistable degeneration,
hence the surjectivity of the map $\Gamma(\alpha)$ is a consequence of the
Clemens local invariant cycle theorem (see \cite{Cl77}). Since by
retraction
$H^{q}(\overline{Pic}^{i}(\mathcal{X}))=H^{q}(\hat{g}^{-1}(0))$,
the injectivity of $\Gamma(\alpha)$ follows if we prove that
\begin{equation}\label{diseguaglianza}
{\rm dim}(H^{q}(\hat{g}^{-1}(0)))\le
{\rm dim}(H^{q}(\hat{g}^{-1}(\frac{1}{2}))^{\pi_{1}(\Delta^{*})}).
\end{equation}
The second term of this inequality can be computed by the
Picard-Lefschetz formula. In fact, denoting by $\delta$ the
vanishing cycle, the Picard-Lefschetz formula says that the
generator $\gamma$ of $\pi_{1}(\Delta^{*})$ acts on
$H_{1}(g^{-1}(\frac{1}{2}))$ sending $\beta$ to
$A_{\gamma}(\beta)=\beta + <\alpha\cap\beta>\alpha$. Since the
action of $\gamma$ on $(H^{q}(\hat{g}^{-1}(\frac{1}{2}))$ is given
by $\wedge^{q}A_{\gamma}^{\vee}$ and the matrix of $A_{\gamma}$
in a suitable basis is given by $$M=\left(\begin{matrix}
1 & 1 & 0 & 0 & \ldots & 0 \\
0 & 1 & 0 & 0 & \ldots & 0 \\
0 & 0 & 1 & 0 & \ldots & 0 \\
0 & 0 & 0 & \ldots & \ldots & 0 \\
0 & 0 & 0 & 0 & \ldots & 0 \\
0 & 0 & 0 & 0 & \ldots & 1 \\
\end{matrix}\right),$$
we have
${\rm dim}(H^{1}(\hat{g}^{-1}(\frac{1}{2}))^{\pi_{1}(\Delta^{*})})=2p-1$
where $p$ is the genus of the general fiber of $g$. A
straightforward computation   also yields
${\rm dim}(H^{2}(\hat{g}^{-1}(\frac{1}{2}))^{\pi_{1}(\Delta^{*})})
=2p^{2}-3p+2$.

The first term of (\ref{diseguaglianza}) is determined by using the
known description of the compactified Picard group of a curve with
a node. The variety $\hat{g}^{-1}(0)$ is obtained starting from a
$\mathbb{P}^{1}$-bundle  over the Jacobian $J$ of the
normalization  of $g^{-1}(0)$. This $\mathbb{P}^{1}$-bundle has
two preferred sections and $\hat{g}^{-1}(0)$ is obtained
identifying the two sections by a translation on $J$. It follows
that $\hat{g}^{-1}(0)$ is homeomorphic to a topologically locally
trivial bundle over $J$ whose fiber $F$ is obtained from a sphere
by identifying two points. By the Leray spectral sequence we deduce
$${\rm dim}(H^{1}(\hat{g}^{-1}(0)))\le \sum_{i+j=1}{\rm dim}(H^{i}(F)\otimes
H^{1}(j))=2p-1$$
$${\rm dim}(H^{2}(\hat{g}^{-1}(0)))\le \sum_{i+j=2}{\rm dim}(H^{i}(F)\otimes
H^{1}(j))=2p^{2}-3p+2.$$
\end{proof}

Since locally over small neighborhoods of points in
$U^{0}\setminus U^{s}$ the families $\Psi$ and $\Phi$ are
homeomorphic to families of the form $\hat{g}\times
id:\overline{Pic}^{i}(\mathcal{X})\times\Delta^{4}\rightarrow
\Delta^{5}$ where
$\hat{g}:\overline{Pic}^{i}(\mathcal{X})\rightarrow \Delta$ is as
in the previous lemma  and $id:\Delta^{4}\rightarrow \Delta^{4}$
is the identity, statement 2 follows from Lemma \ref{lemma
degenerazione}.

It remains to prove statement b).
Its proof is a slight modification of the proof of the degeneration of
the Leray spectral sequence of a smooth projective fibration.
%Now we prove that, for $p+q=2$, the $E_{2}^{p,q}$
%terms of the Leray spectral sequences associated with $\overline{\Phi}$
%and $\overline{\Psi}$ survive to the $E_{\infty}$ page.
Since the same
proof works for both $\Psi_{0}$ and $\Phi_{0}$, we deal
explicitely only with the
first case.

Let $d_{l}^{p,q}:E_{l}^{p,q}\rightarrow E_{l}^{p+l,q-l+1}$ be the
differentials on the l-th page of the spectral sequence. We need
to prove that $d_{3}^{0,2}=0$ and $d_{2}^{p,q}=0$ for $p+q\le
2$. Obviously $d_{2}^{p,0}=0$.

Now we prove  $d_{2}^{p,1}=0$. For any
$\alpha\in H^{2}(\widetilde{M}_{(0,2,2)},\Q)$
%be the class of an ample
%divisor. The
the cup product with the restrictions of
$\wedge^{j}\alpha$ to the fibers of $\Psi_{0}$ induces  morphisms
$(\wedge^{j}\alpha)^{q}:R^{q}\Psi_{0*}\Q\rightarrow
R^{q+2j}\Psi_{0*}\Q$.
If $\alpha$ is the class of an ample divisor
the irreducibility of the fibers of $\Psi_{0}$
implies that
%Since the fibers of $\Psi_{0}$ are
%irreducible,
the morphism $(\wedge^{5}\alpha)^{0}$ is an
isomorphism. On the other hand, for $q>0$, the morphism
$(\wedge^{5}\alpha)^{q}:R^{q}\Psi_{0*}\Q\rightarrow
R^{q+10}\Psi_{0*}\Q$ is trivially $0$. By Lemma 4.13 of
\cite{Vo03} the differentials $d_{l}^{p,q}$ commute with the maps
$H^{p}((\wedge^{j}\alpha)^{q}):H^{p}(R^{q}\Psi_{0*}\Q)\rightarrow
H^{p}(R^{q+2j}\Psi_{0*}\Q)$ induced in cohomology by the maps
$(\wedge^{j}\alpha)^{q}$, therefore we get
$H^{p+2}((\wedge^{5}\alpha)^{0})\circ d_{2}^{p,1}=
d_{2}^{p,11}\circ H^{p}((\wedge^{5}\alpha)^{1})$. It follows
$d_{2}^{p,1}=0$.

It remains to show that $d_{2}^{0,2}=0$ and $d_{3}^{0,2}=0$.
%We
%notice that $H^{0}(R^{2}\Psi_{0*}\Q)$ has dimension one.
%In fact,
By the global invariant cycle theorem (see Theorem 4.24 of
\cite{Vo03}) any section in $ H^{0}(R^{2}\Psi_{s*}\Q)$
is the image of a global cohomology class $\alpha\in H^{2}(\Mdt,\Q)$:
more precisely it is in the image  of
$H^{0}(i^{*})\circ H^{0}((\wedge^{1}\alpha)^{0})$.
By statement 2 the map
$H^{0}(i^{*}):H^{0}(R^{2}\Psi_{0*}\Q)\rightarrow
%H^{0}(R^{2}\Psi_{s*}\Q)\simeq
H^{0}(i^{*}R^{2}\Psi_{0*}\Q)\simeq
H^{0}(R^{2}\Psi_{s*}\Q)$ is an isomorphism:
hence any section in  $H^{0}(R^{2}\Psi_{0*}\Q)$
is in the image of  $H^{0}((\wedge^{1}\alpha)^{0})$ for a certain
$\alpha\in  H^{2}(\widetilde{M}_{(0,2,2)},\Q)$.
%Letting $x\in U^{s}$, the vector space
%$H^{0}(R^{2}\Psi_{s*}\Q)$ is isomorphic to the subspace
%$H^{2}(\Psi^{-1}(x),\Q)^{\pi_{1}(U^{0})}\subset
%H^{2}(\Psi^{-1}(x),\Q)$ invariant under the monodromy action of
%$\pi_{1}(U^{0})$.
%By the global invariant cycle theorem (see Theorem 4.24 of
%\cite{Vo03}) $H^{2}(\Psi^{-1}(x),\Q)^{\pi_{1}(U^{0})}$
%equals the
%image of $H^{2}(\mathbf{M},\Q)$ in $H^{2}(\Psi^{-1}(x),\Q)$. Since
%$\Psi^{-1}(x)$ is a Lagrangian subvariety, this image is a sub
%Hodge structure included in $H^{1,1}_{\Q}(\Psi^{-1}(x))$.
%Since
%any hyperelliptic curve of genus 5 can be obtained as a curve in a
%linear system of the form $|2H|$ for $H$ a degree $2$ ample divisor
%on a K3 surface, we can suppose $\Psi^{-1}(x)$ isomorphic to the Jacobian
%of a general hyperelliptic curve: hence
%$H^{1,1}_{\Q}(\Psi^{-1}(x))\simeq\Q$ (see \cite{Pi88}).
%It follows that the vector space $H^{0}(R^{2}\Psi_{0*}\Q)$ is generated by the
%image of the morphism $H^{0}((\wedge^{1}\alpha)^{0})$.
Using again Lemma 4.13 of \cite{Vo03} we get the equality
$H^{2}((\wedge^{1}\alpha)^{-1})\circ d_{2}^{0,0}= d_{2}^{0,2}\circ
H^{0}((\wedge^{1}\alpha)^{0})$
and, since $d_{2}^{0,0}=0$, we obtain $d_{2}^{0,2}=0$.
An analogous argument proves also $d_{3}^{0,3}=0$.\qed\medskip

\noindent{\em Proof of Proposition \ref{coomologia aperto}.}
%\begin{proof}
The Proposition is a consequence of the existence of an open smooth
irreducible dense subset $V\subset \Phi^{-1}(R)$ such that
$H^{1}(V,\Q)=0$.
In fact, assuming the existence of such a $V$, the couple
$(\Phi^{-1}(U)\cup V,\Phi^{-1}(U))$ induces the exact sequence
$$H^{2}(\Phi^{-1}(U)\cup V,\Phi^{-1}(U))\stackrel{a}{\rightarrow}
H^{2}(\Phi^{-1}(U)\cup V)\rightarrow
H^{2}(\Phi^{-1}(U))\cdots$$
$$\cdots\rightarrow H^{3}(\Phi^{-1}(U)\cup
V,\Phi^{-1}(U)).$$ Since $V$ is smooth, by excision and Thom
isomorphism we get $H^{i}(\Phi^{-1}(U)\cup
V,\Phi^{-1}(U))=H^{i-2}(V)$ for $i\ge 2$. Hence
$H^{3}(\Phi^{-1}(U)\cup V,\Phi^{-1}(U))=0$ and
$H^{i}(\Phi^{-1}(U)\cup V,\Phi^{-1}(U))=\Q$ by irreducibility of
$V$. Since $V\subset \Phi^{-1}(R)$ is dense, the complement of
$\Phi^{-1}(U)\cup V$ in $M_{(0,2,1)}$ has codimension two, therefore
$H^{2}(\Phi^{-1}(U)\cup V)= H^{2}(M_{(0,2,9)})=\Q^{23}$.
Finally $a$ is injective
because its image  contains the Chern class of the line bundle
associated with the effective divisor $\Phi^{-1}(R)$: hence
$H^{2}(\Phi^{-1}(U))=\Q^{22}.$

In order to define $V$ we denote by
$R^{0}\subset R$ the locus parametrizing
curves of the form
$C=C_{1}\cup C_{2}$, where $C_{1}\ne C_{2}$, the singular locus of $C_{i}$
consists of at most a nodal point and
$C_{1}\cap C_{2}$ is included in the smooth locus of both $C_{1}$ and
$C_{2}$.

A concrete description of the inclusion $R^{0}\subset R$ is given by
means of the map
 $f:X\rightarrow \mathbb{P}^{2}$ (see Remark \ref{rivestimento doppio}).
This map
identifies $R^{0}$ with the set of pairs of distinct lines in
$\mathbb{P}^{2}$ such that the intersection of each line with the
branch locus $S$ of $f$ is either reduced or contains at most a unique
double point and, in this second case, the support of this double
point does not belong to the intersection of the two lines.

The subvariety $V\subset \Phi^{-1}(R^{0})$ is defined as the locus
parametrizing sheaves of the form $F=i_{*}L$, where $i:C\rightarrow X$
is the inclusion of a curve of $R^{0}$ and $L$ is a line bundle on $C$.

Before giving a global description of $V$, we study the locus of
$M_{(0,2,1)}$ parametrizing sheaves supported on a fixed curve in
$R^{0}$.
\begin{lem}\label{lemma stabili}
Let $C=C_{1}\cup C_{2}$ be a curve in $R^{0}$. Denote by
$i:C\rightarrow X$  the inclusion of $C$ and by
$i_{1}:C_{1}\rightarrow X$ and $i_{2}:C_{2}\rightarrow X$ the
inclusion of its components. Let $F$ be a torsion free sheaf on
$C$ and suppose $F=i_{*}G\in \Phi^{-1}(C)$.
Then:\begin{enumerate}[1)]\item{Up to exchange of $C_{1}$ and $C_{2}$, the
sheaf $F$ fits in an exact sequence of the form
\begin{equation}\label{succesatta}
0\rightarrow i_{1*}L_{1} \rightarrow F \rightarrow
i_{2*}L_{2}\rightarrow 0
\end{equation}
where $L_{1}$ and $L_{2}$ are rank $1$ torsion free sheaves
whose  degrees are one and two respectively,}
\item{For any non trivial extension of the form (\ref{succesatta})
the middle term $F$ is a stable sheaf,}
\item{Fixing $L_{1}$ and $L_{2}$, two non trivial extensions of the form (\ref{succesatta})
have middle terms isomorphic if and only if they differ by a scalar
multiplication,}
\item{Fixing $L_{1}$ and $L_{2}$, for any point $p$
in $C_{1}\cap C_{2}$ there exists a unique, up
to $\C^{*}$, non trivial extension of the form (\ref{succesatta})
such that the restriction $G$ of $F$ to $C$
 is not locally free at $p$.}
\end{enumerate}
\end{lem}
\begin{proof}
1) Let $G_{1}$ and $G_{2}$ be the torsion free parts of the
restrictions of $G$ to $C_{1}$ and $C_{2}$, then $F$ fits in an exact
sequence of the form
$$0\rightarrow F\rightarrow i_{1*}G_{1}\oplus i_{2*}G_{2}\rightarrow Q
\rightarrow 0 $$ where $Q$ is a quotient of the schematic
intersection of $C_{1}$ and $C_{2}$. Stability and $ch_{2}(F)=1$
imply that either ${\rm deg}(G_{1})={\rm deg}(G_{2})=2$ and ${\rm
length}(Q)=1$ or $\{{\rm deg}(G_{1}),{\rm deg}(G_{2})\}=\{2,3\}$
and ${\rm length}(Q)=2$. Supposing ${\rm deg}(G_{2})=2$ and
setting $G_{2}=L_{2}$ and $L_{1}:=Ker(F\rightarrow i_{2*}L_{2})$
we get the sequence
(\ref{succesatta}).
2) If $F$ were unstable, there would be a sheaf of the form
$i_{j*}M$ with ${\rm deg}(M)=2$ injecting into $F$. If $j=1$ this
would imply that $M$ is a subsheaf of $L_{1}$: absurd. If $j=2$
then $M=L_{2}$ and the sequence splits.
3) By 2) $End(F)=\C$.
Since $Hom(F,i_{2*}L_{2})=End(i_{1*}L_{1})=\C$ a diagram chase
proves 3).
4) By the Grothendieck spectral sequence of $H^{0}\circ
\mathcal{H}om$, since $\mathcal{H}om(i_{2*}L_{2},i_{1*}L_{1})=0$,
we get $Ext^{1}(i_{2*}L_{2},i_{1*}L_{1})=
H^{0}(\mathcal{E}xt^{1}(i_{2*}L_{2},i_{1*}L_{1}))$. The sheaf
$\mathcal{E}xt^{1}(i_{2*}L_{2},i_{1*}L_{1}))$ is isomorphic to the
structure sheaf of the schematic intersection of $C_{1}$ and
$C_{2}$ and the extensions of the form (\ref{succesatta}) with $G$
locally free near $C_{1}\cap C_{2}$ correspond to sections
generating $\mathcal{E}xt^{1}(i_{2*}L_{2},i_{1*}L_{1})$: item 4)
follows.
\end{proof}
It remains to prove that $V$ is  an open, dense, smooth, irreducible
subset of $\Phi^{-1}(R)$ and $H^{1}(V,\Q)=0$.
Openness is obvious. By the previous lemma $V$ is dense in
$\Phi^{-1}(R^{0})$ and, since $R^{0}$ is dense in $R$ and the fibers of
$\Phi$ are equidimensional (see \cite{Ma00}), the open subvariety  $V$
is dense in the
divisor $\Phi^{-1}(R)$.

Smoothness holds because $R^{0}$ is smooth and
the differential of $\Phi$ at any point $p=F$ of $V$ is surjective: in fact
since the restriction of $F$ to its support $C$ is a line bundle, such a
differential is identified with the natural map
$d:Ext^{1}(F,F)\rightarrow H^{0}(\mathcal{E}xt^{1}(F,F))\simeq
H^{0}(N_{C|X})$ and this map is surjective because its cokernel is
always included in $H^{2}(\mathcal{H}om^{1}(F,F))$ which is zero
since $F$ is supported on a curve.

In order to prove the irreducibility of $V$ we show that there
exists a $\mathbb{P}^{1}$-bundle over an irreducible base $b:
P\rightarrow N$ having a surjective map, actually birational,
to $\Phi^{-1}(R^{0})$.
More precisely, denoting by $g:M_{(0,1,0)}\times
M_{(0,1,1)}\rightarrow |H|^{2}$ the product of the maps
$\Phi_{(0,1,0)}$ and $\Phi_{(0,1,1)}$,
defined at the beginning of Section 1, and denoting
by $T\subset |H|^{2}$ the inverse image of
$R^{0}\subset R =Sym^{2}(|H|)$ in $|H|^{2}$, the base $N$ is
$g^{-1}(T)$ and, keeping notation as in Lemma \ref{lemma stabili},
the fiber $b^{-1}(i_{1*}L_{1},i_{2*}L_{2})$ is naturally
isomorphic to $\mathbb{P}(Ext^{1}(i_{2*}L_{2},i_{1*}L_{1}))$.

Any point $p\in M_{(0,1,0)}\times M_{(0,1,1)}$ has an open
neighborhood $U_{p}$ in the classical topology of the form
$U_{1}\times U_{2}$ such that each $X\times U_{i}$ is endowed with
a tautological family $F_{i}$. Let $q_{i}:X\times U_{1}\times
U_{2}\rightarrow X\times U_{i}$  and $q:X\times U_{1}\times
U_{2}\rightarrow U_{1}\times U_{2}$ be the projections. For any
$U_{p}\subset g^{-1}(T)$ the sheaf
$\mathcal{E}xt^{1}_{q}(q_{2}^{*}F_{2},q_{1}^{*}F_{1})$ is a rank
$2$ vector bundle. By 3) of Lemma \ref{lemma stabili} the fibers
of the associated  projective bundle
$b_{U_{p}}:\mathbb{P}(\mathcal{E}xt^{1}_{q}(q_{2}^{*}F_{2},q_{1}^{*}F_{1}))\rightarrow
U_{p}$ parametrize isomorphism classes of sheaves $F$ fitting in a
non trivial extension of the form (\ref{succesatta}): it follows
that the bundles $b_{U_{p}}$ can be glued to form a global
$\mathbb{P}^{1}$-bundle $b: P\rightarrow N$. By 1) and 2) of Lemma
\ref{lemma stabili} the natural modular map $\phi: P\rightarrow
M_{(0,2,1)}$ surjects onto $\Phi^{-1}(R^{0})$: hence the open
subset $V\subset \Phi^{-1}(R^{0})$ is irreducible.

Let $N^{0}\subset N\subset M_{(0,1,0)}\times M_{(0,1,1)}$ be the
open subset parametrizing pairs of sheaves whose restrictions to
their supports are line bundles, set $P^{0}:=b^{-1}(N^{0})$ and
denote by $b^{0}: P^{0}\rightarrow N^{0}$ and by $\phi^{0}:
P^{0}\rightarrow M_{(0,2,1)}$ the restrictions of $b$ and $\phi$.
By 4) of Lemma \ref{lemma stabili} the locus  $W\subset P^{0}$,
parametrizing extensions of the form (\ref{succesatta}) whose middle
terms have locally free restrictions to their supports, is the
complement of a two section $D$.

The map $\phi^{0}$ induces a bijection, hence an isomorphism, between
the smooth varieties $W$ and $V$. In fact,
by 1) of Lemma \ref{lemma stabili},
%implies $\phi^{0}(W)=V$ and, for  $F=i_{*}G\in V$
if $F=i_{*}G\in V$, then $F$
%$\in \phi^{0}(V)$ by 1)
%of Lemma \ref{lemma stabili}.
 is the middle term in an exact
sequence of the form (\ref{succesatta}): hence $\phi^{0}(W)=V$.
Moreover since $G$ is a line
bundle of degree $5$, its restrictions to the components of $C$
have degree $2$ and $3$: hence the sheaf $i_{2*}L_{2}$ in the sequence
(\ref{succesatta}) is
unique quotient of $F$ belonging to $M_{(0,1,1)}$
and $i_{1*}L_{1}\in M_{(0,1,0)} $ is the associated kernel: it follows
that $(\phi^{0})^{-1}(F)\subset (b^{0})^{-1}(i_{1*}L_{1},i_{2*}L_{2})$.
Since, by 3) of Lemma \ref{lemma
stabili} the map $\phi^{0}$ is injective on the fibers of $b^{0}$,
the restriction of  $\phi^{0}$ to $W$ is injective too.

% bijection (hence an isomorphism)
%between the smooth varieties $W$ and $V$.

We now show that $H^{1}(W)=H^{1}(V)$ is zero.
Let $D_{s}$ be the smooth locus of $D$.
The couple $(W\cup D_{s},W)$ induces the long exact sequence
$$H^{1}(W\cup D_{s})\rightarrow H^{1}(W)\rightarrow H^{2}(W\cup D_{s},W)
\stackrel{a}{\rightarrow} H^{2}(W\cup D_{s}).$$ The vector space
$H^{1}(W\cup D_{s})$ is zero. In fact the complement of $W\cup
D_{s}$ in $P^{0}$ has codimension two, hence $H^{1}(W\cup
D_{s})=H^{1}(P^{0})$ and since $P^{0}$ is a
$\mathbb{P}^{1}$-bundle over $N^{0}$ its first cohomology group is
trivial if $H^{1}(N^{0})=0$. The last equality holds since the
complement of $N^{0}$ in the simply connected manifold
$M_{(0,1,0)}\times M_{(0,1,1)}$ has codimension two. Indeed it is
the union of $g^{-1}(|H|^{2}\setminus T)$ and the locus $Y\subset
M_{(0,1,0)}\times M_{(0,1,1)}$ parametrizing pairs of sheaves of
the form $(i_{1*}L_{1},i_{2*}L_{2})$ where either $L_{1}$ or
$L_{2}$ is not a line bundle. The subvariety
$g^{-1}(|H|^{2}\setminus T) \subset M_{(0,1,0)}\times M_{(0,1,1)}$
has codimension two since $|H|^{2}\setminus T$ has codimension two
in $|H|^{2}$ and by \cite{Ma00} the fibers of $g$ are
equidimensional. The locus $Y\subset M_{(0,1,0)}\times
M_{(0,1,1)}$ has codimension two since the same property holds for
the subsets of $M_{(0,1,0)}$ and $M_{(0,1,1)}$ parametrizing
sheaves which are not push-forwards of line bundles from curves in
$|H|$.

Using the previous exact sequence, it remains to show the
injectivity of $a$. The image of $a$ is the Chern class of the
line bundle associated with the divisor $D_{s}$ and it is not zero
since $D$ has degree $2$ on the fibers of $b^{0}$. Since by
excision and Thom isomorphism the dimension of $H^{2}(W\cup
D_{s},W)$ is the number of connected components of $D_{s}$ we need
to prove that $D_{s}$ is connected or, equivalently, $D$ is
irreducible.

Denote by $Z\subset X\times T\subset X\times |H|^{2}$ the
incidence subvariety parametrizing triplets of the form
$(p,C_{1},C_{2})$ where $p\in C_{1}\cap C_{2}$.
There exists a regular morphism
$m:D\rightarrow Z$ given by sending
 a
non trivial extension of the form (\ref{succesatta}), where
$F=i_{*}G$ and  $G$ is not locally free, to the triplet
$(p,C_{1},C_{2})$ where $p\in C_{1}\cap C_{2}$ is the unique point
at which $G$ is not locally free. By 4) of Lemma \ref{lemma stabili},
the fiber of $m$  over
$(p,C_{1},C_{2})$  is isomorphic to  $Pic^{1}(C_{1})\times
Pic^{2}(C_{2})$ hence it is irreducible and of constant dimension $4$:
therefore the irreducibility of  $D$ follows from the
one of $Z$. Finally $Z$ is irreducible since  the
projection $p:Z\rightarrow |H|^{2}$ is a double covering and it is
obtained from the double covering $f:X\rightarrow \mathbb{P}^{2}$
by base change with the map $q:T\rightarrow \mathbb{P}^{2}$
sending $(C_{1},C_{2})$ to $f(C_{1}\cap C_{2})$: since $X$ and $T$
are irreducible also $Z$
is irreducible \qed\medskip

\noindent{\em Proof of Proposition \ref{componenti irriducibili}.}
Since the fibers of $\Psi$ are equidimensional
and the exceptional divisor of O'Grady's desingularization is
irreducible, we need to prove  the irreducibility of the stable
locus of $\Phi_{(0,2,2)}^{-1}(R^{00})$, where $R^{00}$
is an open dense subset of $R$.
It is irreducible since $R$ is
irreducible and for general  $C=C_{1}\cup C_{2} \in R $,  the
stable locus of  $\Phi_{(0,2,2)}^{-1}(C)$ is  a $\C^{*}$-bundle
over $Pic^{3}(C_{1})\times Pic^{3}(C_{2})$ (this can be proved  as
in the case  $C\in R(1)\cup R(2)$ in Proposition 2.1.4. of
\cite{Ra04}).\qed

\section{A basis for $H^{2}(\mathbf{M}, \Z)$}
Following \cite{OG99} we denote by $\St\subset\M$ the exceptional
divisor of the desingularization map $\widetilde{\pi}:\M\rightarrow
M_{(2,0,-2)}$ and by $\widetilde{B}\subset\M$ the strict transform of
the irreducible divisor $B\subset M_{(2,0,-2)}$ parametrizing not
locally free sheaves. Let $M^{U}$ be the Uhlenbeck compactification of
the $\mu$-stable locus of $M_{(2,0,-2)}$ and let
$\varphi:M_{(2,0,-2)}\rightarrow M^{U}$ be the functorial morphism. Let
$\mu:H^{2}(X,\Z)\rightarrow H^{2}(M^{U},\Z)$ be the Donaldson morphism
(see \cite{Li93},\cite{Mo93} and \cite{FM94})
and set $\widetilde{\mu}:=\widetilde{\pi}^{*}\circ\varphi^{*}\circ\mu$.
\begin{thm}\label{base intera}
The morphism $\widetilde{\mu}$ is injective and
$$H^{2}(\M,\Z)=\widetilde{\mu}(H^{2}(X,\Z))\oplus\Z c_{1}(\St)\oplus\Z c_{1}(\widetilde{B}).$$
\end{thm}
\begin{proof}
By $\S$ 5 of \cite{OG99}   $\widetilde{\mu}$ is
injective and
$\widetilde{\mu}(H^{2}(X,\Z))\oplus\Z c_{1}(\St)\oplus\Z c_{1}(\widetilde{B}) \subset
H^{2}(\M,\Z)$ .
We first prove that the submodule
$\widetilde{\mu}(H^{2}(X,\Z))$ is saturated.
In fact by \cite{FM94}, VII.2.17
(see also formula 5.1 of \cite{OG99})
\begin{equation}\label{formula}
\int_{\M}\bigwedge^{10}\widetilde{\mu}(\alpha)=\frac{10!}{5!2^{5}}(\int_{X}\bigwedge^{2}\alpha )^{5}.
\end{equation}
Polarizing this formula we get
\begin{equation}\label{formula polarizzata}
\int_{\M}\bigwedge_{i=1}^{10}\widetilde{\mu}(\alpha_{i})=\frac{10!}{5!2^{5}}\frac{1}{10!}\sum_{\sigma\in S_{10}}
(\Pi_{i=1}^{5}\int_{X}\alpha_{\sigma(2i-1)}\wedge\alpha_{\sigma(2i)}).
\end{equation}
Let $\alpha_{1}\in H^{2}(X,\Z)$ be a primitive element, since
$(H^{2}(X,\Z),\wedge)$ is unimodular of dimension 24, there exist
$\alpha_{2},\ldots,\alpha_{10}\in H^{2}(X,\Z)$ such that
\[\begin{array}{ccc}
\int_{X}\alpha_{j}\wedge\alpha_{k}=1 & if & \{j,k\}=\{2i-1,2i\} \\
\int_{X}\alpha_{j}\wedge\alpha_{k}=0 & if & \{j,k\}\ne\{2i-1,2i\}.
\end{array}\]
By Formula (\ref{formula polarizzata}) we have
$\int_{\M}\bigwedge_{i=1}^{10}\widetilde{\mu}(\alpha_{i})=1$. Hence
$\widetilde\mu$ sends primitive elements of $H^{2}(X,\Z)$ to primitive elements of
 $H^{2}(\M,\Z)$, therefore
$\widetilde\mu(H^{2}(X,\Z))$ is saturated.

We now prove that $\widetilde{\mu}(H^{2}(X,\Z))\oplus\Z
c_{1}(\St)\oplus\Z c_{1}(\widetilde{B})$ is saturated by
evaluating a basis on suitable homology classes. By Proposition
3.0.5 of \cite{OG99}, there exists an open dense subset
$B^{0}\subset B$ which is a $\mathbb{P}^{1}$-bundle over the
smooth locus of the symmetric product $X^{(4)}$ and by the proof
of Lemma 3.0.13 of \cite{OG99} the intersection $\Sigma\cap B^{0}$
is a three-section of this $\mathbb{P}^{1}$-bundle. Since by
\cite{LS05} the smooth variety $\M$ is the blow up of
$M_{(2,0,-2)}$ along $\Sigma$,  the strict transform
$\widetilde{B}$ has an open subset which is a
$\mathbb{P}^{1}$-bundle and, denoting by $\gamma$ its fiber, we
have $\gamma\cdot\widetilde{\Sigma}=3$. On the other hand $\St$
has an open dense subset $\St^{0}$ which is
$\mathbb{P}^{1}$-bundle on the smooth locus of $\Sigma$ and,
denoting by $\delta$ a fiber, we get $\delta\cdot\widetilde{B}=1$.
Finally, since $\M$ has trivial canonical bundle, we also get
$\gamma\cdot\St=-2=\delta\cdot\widetilde{B}$.

Let $(\alpha_{1},\ldots,\alpha_{22})$ be a basis of $H^{2}(X,\Z)$,
since $\widetilde{\mu}(H^{2}(X,\Z))$ is saturated there exist
$\beta_{1},\ldots,\beta_{22}\in H_{2}(\M,\Z)$ such that ${\rm
det}(\alpha_{i}(\beta_{j}))_{i,j\le 22}=1$. Moreover
$\alpha_{i}(\delta)=\alpha_{i}(\gamma)=0$ for any i, since
$\delta$ and $\gamma$ are contracted by
$\varphi\circ\widetilde{\pi}$ (see Proposition (3.0.5) of \cite{OG99}).
Therefore, denoting by $M$ the
evaluation matrix of
$(\alpha_{1},\ldots,\alpha_{22},c_{1}(\St),c_{1}(\widetilde{B}))$
on $(\beta_{1},\ldots\beta_{22},[\delta],[\gamma])$, we have
$${\rm
det}(M)={\rm det}\left(
\begin{array}{cc}
c_{1}(\St)([\delta]) & c_{1}(\widetilde{B})([\delta])  \\
c_{1}(\St)([\gamma]) & c_{1}(\widetilde{B})([\gamma])
\end{array}
\right)={\rm det}\left(
\begin{array}{cc}
-2 & 1  \\
3  & -2
\end{array}
\right)=1$$ which implies the statement.
\end{proof}

\section{The Beauville form of $\M$}

In this section  we determine the Beauville form and the Fujiki
constant of $\M$. Before analyzing the case of $\M$, we recall the
theorem due to Beauville and Fujiki (see \cite{Be83,Fu87}) that
defines the Beauville form and the Fujiki constant of an
irreducible symplectic variety.

\begin{thm}\label{teorema def forma beauville}
Let $Y$ be   irreducible symplectic variety of dimension $2n$.
There exist a unique indivisible bilinear integral symmetric form
$B_{Y}\in S^{2}(H^{2}(Y,\mathbb{Z}))^{*}$, called the Beauville
form, and a unique positive constant $c_{Y}\in\mathbb{Q}$, called
the Fujiki constant, such that for any $\alpha\in
H^{2}(Y,\mathbb{C})$
\begin{equation}\label{formula fujiki}\int_{Y}\alpha^{2n}=c_{Y}B_{Y}(\alpha,\alpha)^{n}\end{equation}
and for $0\neq\omega\in H^{0}(\Omega^{2}_{Y})$
\begin{equation}B_{Y}(\omega+\overline{\omega},\omega+\overline{\omega})>0.
\end{equation}\end{thm}

\begin{rem} Formula (\ref{formula fujiki}) is Fujiki's formula.
We will use also its  polarized form:
\begin{equation}\label{formula fujiki polarizzata}
\int_{Y}\alpha_{1}\wedge...\wedge\alpha_{2n}=\frac{c_{Y}}{2n!}
\sum_{\sigma\in
S_{2n}}B(\alpha_{\sigma(1)},\alpha_{\sigma(2)})...B(\alpha_{\sigma(2n-1)},\alpha_{\sigma(2n)})
.\end{equation}
\end{rem}

In the case of $\M$ we have the following theorem.

\begin{thm}\label{mainthm}
Set
$\Lambda:=\mathbb{Z}(c_{1}(\widetilde{\Sigma}))\oplus\mathbb{Z}c_{1}(\widetilde{B})\subset
H^{2}(\M,\mathbb{Z})$. The direct sum decomposition
\[H^{2}(\M,\mathbb{Z})=
\widetilde{\mu}(H^{2}(X,\mathbb{Z}))\oplus_{\bot} \Lambda\] is
orthogonal with respect to
$B_{\M}$.\\
The  map
$\widetilde{\mu}:(H^{2}(X,\mathbb{Z}),\wedge)\longrightarrow
(H^{2}(\M,\mathbb{Z}),B_{\M})$
is an isometric embedding.\\
The matrix of the Beauville form on $\Lambda$  is
\[
\begin{array}{|c|c|c|}
\hline
&c_{1}(\widetilde{\Sigma})  & c_{1}(\widetilde{B}) \\
\hline
c_{1}(\widetilde{\Sigma}) & -6 & 3 \\
\hline
c_{1}(\widetilde{B})  & 3 & -2 \\
\hline
\end{array}\; .
\]
\end{thm}

In order to determine the Beauville form $B_{\M}$ of $\M$ we use the
following fundamental property of the Donaldson morphism
(see \cite{Li93},\cite{Mo93} and \cite{FM94}).

 Let $Y$ be a smooth variety and let $\mathcal{F}$ be a
coherent sheaf on $X\times Y$. Suppose that $\mathcal{F}$ is a flat
family parametrizing semistable sheaves on $X$ with Mukai vector
$(2,0,-2)$. Let $\iota_{\mathcal{F}}:Y\rightarrow M_{(2,0,-2)}$ be the
associated modular map and let $p$ and $q$ be the projections of
$X\times Y$ on $X$ and $Y$ respectively.
For any $\alpha\in H^{2}(X,\Z)$,
\begin{equation}\label{proprieta Donaldson}
\iota_{\mathcal{F}}^{*}\circ\varphi^{*}\circ\mu(\alpha)=
q_{*}(c_{2}(\mathcal{F})\wedge p^{*}(\alpha)).
\end{equation}
\begin{lem}\label{conti}
Let $\omega$ be a symplectic holomorphic two form on $X$.
\begin{enumerate}
\item{$\int_{\M}\Muo^{8}\wedge c_{1}(\widetilde{B})\wedge
c_{1}(\St)=35\cdot 9(\int_{X}\om^{2})^{4}$}
\item{$\int_{\M}\Muo^{8}\wedge
c_{1}(\St)^{2}=-35\cdot 9\cdot 2(\int_{X}\om^{2})^{4}$}
\item{$\int_{\M}\Muo^{8}\wedge c_{1}(\widetilde{B})^{2}=-35\cdot
3\cdot2(\int_{X}\om^{2})^{4}$}
\end{enumerate}
\end{lem}
\begin{proof}
Let $\mathcal{I}$ be the universal sheaf of ideals on $X\times
X^{[2]}$ and let $\pi_{i}:X\times X^{[2]}\times X^{[2]}\rightarrow
X\times X^{[2]}$ be the projection on the product of the first and
the i-th factors. Denote by $\iota$ the modular map associated
with the family
$\pi_{2}^{*}(\mathcal{I})\oplus\pi_{3}^{*}(\mathcal{I})$. By
Formula (\ref{proprieta Donaldson})
$$\iota^{*}\circ\varphi^{*}\circ\mu(\alpha)=
p_{1}^{*}(\alpha^{[2]})+p_{2}^{*}(\alpha^{[2]})$$ where
$p_{j}:X^{[2]}\times X^{[2]}\rightarrow X^{[2]}$ is the projection
on the j-th factor and $\alpha^{[2]}\in H^{2}(X^{[2]},\C)$ is the
class associated with $\alpha$, namely the pull-back from
$X^{(2)}$ of the class whose pull-back to $X^{2}$ restricts to
$\alpha$ on each component. It follows that $$\int_{X^{[2]}\times
X^{[2]}}(\iota^{*}\circ\varphi^{*}\circ\mu\om)^{8}=\int_{X^{[2]}\times
X^{[2]}}((p_{1}^{*}+p_{2}^{*})(\omega^{[2]}+\overline{\omega}^{[2]}))^{8}=$$
$$70\int_{X^{[2]}\times
X^{[2]}}(p_{1}^{*}(\omega^{[2]}+\overline{\omega}^{[2]}))^{4}\wedge
(p_{2}^{*}(\omega^{[2]}+\overline{\omega}^{[2]}))^{4}=
70(\int_{X^{[2]}}(\omega^{[2]}+\overline{\omega}^{[2]})^{4})^{2}=$$
$$70\cdot 9(\int_{X}\om^{2})^{4}.$$
Since $\iota$ is a double covering over $\Sigma$ and
$\widetilde{B}\cap\St^{0}$ is a
rational section  of the restriction of $\varphi\circ\widetilde{\pi}$ to $\St^{0}$ we have $\int_{\M}\Muo^{8}\wedge c_{1}(\widetilde{B})\wedge
c_{1}(\St)=35\cdot 9(\int_{X}\om^{2})^{4}.$
Since the normal bundle to $\St$ has degree $-2$ on the fibers of
$\St^{0}$, we also get $\int_{\M}\Muo^{8}\wedge
c_{1}(\St)^{2}=-35\cdot 9\cdot 2(\int_{X}\om^{2})^{4}$.
Analogously, since $\widetilde{B}^{0}\cap\St$ is a rational three-section of
the restriction of $\varphi\circ\widetilde{\pi}$ to $\widetilde{B}^{0}$ and the normal bundle of $\widetilde{B}$ has
degree $-2$ on the general fiber, we get
$\int_{\M}\Muo^{8}\wedge c_{1}(\widetilde{B})^{2}=-35\cdot
3\cdot 2(\int_{X}\om^{2})^{4}$.
\end{proof}

\noindent{\em Proof of Theorem \ref{mainthm}}
By equality (\ref{formula}) and Fujiki's formula
there exists $a\in\Q$ such that
$B_{\M}(\widetilde{\mu}(\alpha),\widetilde{\mu}(\alpha))=a\int_{X}\alpha^{2}$ for any $\alpha\in H^{2}(X,\Z)$ and
$c_{\M}=\frac{10!}{5!2^{5}\cdot a^{5}}$. Since $\widetilde{B}$ and
$\St$ are contracted by $\varphi\circ\widetilde{\pi}$, we get
$\int_{\M} c_{1}(\widetilde{B}) \wedge\widetilde{\mu}(\alpha)^{9}=
\int_{\M} c_{1}(\St) \wedge\widetilde{\mu}(\alpha)^{9}=0$, hence
$B_{\M}(\widetilde{\mu}(\alpha), c_{1}(\widetilde{B}))=B_{\M}(\widetilde{\mu}(\alpha),c_{1}(\St))=0$.
Using the polarized form of Fujiki's Formula and item 1 of
Lemma \ref{conti}, we get $$\int_{\M}\Muo^{8}\wedge c_{1}(\widetilde{B})\wedge
c_{1}(\St)=$$ $$\frac{8!5\cdot 2\cdot
c_{\M}}{10!}B_{\M}(c_{1}(\widetilde{B}),c_{1}(\St))B_{\M}(\Muo,\Muo)^{4}=
35\cdot
9(\int_{X}\om^{2})^{4}.$$ Hence, comparing the second and the third
term of this equation and replacing $c_{\M}$ and $ B_{\M}(\Muo,\Muo)^{4}$ by their
values,
$$B_{\M}(c_{1}(\widetilde{B}),c_{1}(\St))=3a.$$
Analogously we get
$$B_{\M}(c_{1}(\St),c_{1}(\St))=-6a,$$
$$B_{\M}(c_{1}(\widetilde{B}),c_{1}(\widetilde{B}))=-2a.$$
Since $B_{\M}$ is integral and primitive and $B_{\M}(\Muo,\Muo)>0$
we conclude $a=1$ and $c_{\M}=\frac{10!}{5!2^{5}}$.
\qed

\end{document}